\documentclass[letterpaper, 10 pt, conference]{ieeeconf}  
\IEEEoverridecommandlockouts                              

\let\labelindent\relax
\usepackage{amsmath} 
\usepackage{amssymb}  

\usepackage{tabulary}
\usepackage[english]{babel}
\usepackage{blindtext}
\usepackage[font=footnotesize]{caption}
\usepackage[makeroom]{cancel}
\usepackage{amsfonts}
\usepackage{graphicx}
\usepackage{pgf,tikz}
\usepackage{tabularx}
\usepackage{float}
\usepackage{cite}
\usepackage{enumitem} 
\usepackage{color}
\usepackage{graphicx}
\usepackage{amsfonts}
\usepackage{dsfont}
\usepackage{tikz}
\usepackage{caption}
\usepackage{subcaption}
\usepackage{amsthm}
\usepackage{nicematrix}
\usepackage{algorithm2e}
\RestyleAlgo{ruled}
\usepackage{algpseudocode}
\usepackage[strict]{changepage}

\newcommand{\vect}[1]{\boldsymbol{#1}}

\newcommand{\mc}{\mathcal}
\newcommand{\ov}{\overline}
\newcommand{\ds}{\displaystyle}
\newcommand{\R}{\mathbb{R}}
\newcommand{\Var}{\textrm{Var}}

\newcommand{\be}{\begin{equation}}
	\newcommand{\ee}{\end{equation}}
\newcommand{\ba}{\begin{array}}
	\newcommand{\ea}{\end{array}}
\allowdisplaybreaks

\newtheorem{remark}{Remark}
\newtheorem{theorem}{Theorem}
\newtheorem{lemma}{Lemma}

\newtheorem{definition}{Definition}
\newtheorem{corollary}{Corollary}
\newtheorem{proposition}{Proposition}

\newtheorem{example}{Example}

\title{Optimal selection of the most informative nodes \\for a noisy DeGroot model with stubborn agents
}

\author{Roberta Raineri, Giacomo Como, Fabio Fagnani
\thanks{R. Raineri is with the Department of Electronics and Telecommunications, Politecnico di Torino, Torino, Italy (\texttt{roberta.raineri@polito.it}). G. Como and F. Fagnani are with the Department of Mathematical Sciences, Politecnico di Torino, Torino, Italy (\texttt{\{giacomo.como,fabio.fagnani\}@ polito.it}).  	G. Como is also with the Department of Automatic Control, Lund University, Lund, Sweden.  }
}

\begin{document}

	\maketitle
	\thispagestyle{empty}

	\begin{abstract}
Finding the optimal subset of individuals to observe in order to obtain the best estimate of the average opinion of a society is a crucial problem in a wide range of applications, including policy-making, strategic business decisions, and the analysis of sociological trends. We consider the opinion vector $X$ to be updated according to a DeGroot opinion dynamical model with stubborn agents, subject to perturbations from external random noise, which can be interpreted as transmission errors. The objective function of the optimization problem is the variance reduction achieved by observing the equilibrium opinions of a subset $\mathcal{K}\!\subseteq\!\mathcal{V}$ of agents. We demonstrate that, under this specific setting, the objective function exhibits the property of submodularity. This allows us to effectively design a Greedy Algorithm to solve the problem, significantly reducing its computational complexity. Simple examples are provided to validate our results.
	\end{abstract}
	
	\section{Introduction}
Accurately estimating the average opinion of a society is a crucial challenge, given its implications for policy-making, strategic business decisions and study of sociological trends. Seminal works, such as \cite{Page_Shapiro_1983}, show how public opinion plays a critical role in creating policies that guide decisions and actions in areas such as government, business or society.  Moreover, accurate forecasting may also be helpful in trends and misinformation spread monitoring.  
Typically, opinion polls are a valuable tool commonly used to effectively capture public opinion. However, since it is not feasible to survey the entire population, the main challenge lies in selecting a group of individuals that best represents the population of interest. This problem can be considered more in general as an instance of a subset selection problem. This is a key question in multiple data-driven applications, spanning a wide range of problems, including for instance feature selection in machine learning \cite{Kempe}\cite{KempeApproxSubmod} \cite{5693974},  sensor placement for environmental monitoring \cite{Krause2}, or even smart testing problem, which seeks to develop effective strategies for controlling the spread of epidemics \cite{smarttesting}.

Properly formalizing the problem, it involves choosing $k<n$ variables $X_i$ out of a large set  $n$ in order to best infer a variable $Y$, more precisely corresponding in such case to the average of $X$.
Given prior knowledge of the covariances among $X_i$ and $Y$, the goal is to select the subset that optimally predicts $Y$. As already stated in our preliminary work \cite{IFAC},  a natural way to formalize the optimization problem is to consider as target function the variance reduction on $Y$, namely the difference between the variance of $Y$ and the variance of $Y$ conditioned on the observed variables, and let the selected optimal subset be the one that maximizes this quantity. We notice how this problem is known to be NP-complete \cite{Kempe} and not submodular in general \cite{IFAC}. 

In this work, we focus on a particular instance of the problem {first introduced in \cite{IFAC}},  in the context of a social network where the observable variables are the opinions of the individuals and the quantity to be estimated is the average opinion $Y$ of the community. The individuals $\mc V$ of the society may in particular be distinguished  as stubborn individuals, that keep their opinion unchanged, and regular ones, which are open to social interactions \cite{9539139}.
The opinions of stubborn agents are here assumed to be known and can be referred for instance to any influencial agent, like politicians or authorities. The opinions of regular individuals are assumed to be determined by network interactions plus an additional random noise, according to a DeGroot model \cite{DeGroot74} with stubborn agents \cite{Acemoglu.ea:2013,Como.Fagnani:2016}. Referring to \cite{Huang2013_noisyInformation}, \cite{Banerjee2024_noiseinterpretation} and \cite{Jackson2022_noiseinterpretation}, we outline multiple real-world scenarios that could explain the role of the noise term. It may arise from transmission errors, e.g. an agent receiving the information via a third party acting as a distorting trasmission channel, such as newspapers or TV news reports. Alternatively, it could result from ambiguity introduced by stubborn agent themselves, e.g. a politician deliberately exposing a deceiving version of their opinion during a public speech.

This case study reveals to be significant not only from an interpretative perspective, but also since the variance reduction function exhibits submodularity under this specific model design without requiring any additional assumption on the covariance distribution of random noise. The submodularity property is crucial because, as proven by \cite{Nemhauser}, it guarantees a near-optimal solution adopting a greedy approximation for maximizing the function of interest subject to a cardinality constraint. 

Thus, our contribution is threefold. At first, we extend the preliminary results presented in \cite{IFAC} to the new setting deriving an explicit formulation for the variance reduction function, which acts as the objective in the optimization problem. Next, we prove that this function is submodular in the context of interest, enabling the use of a greedy approximation. Finally, we formalize and apply the Greedy Algorithm to two simple case studies, supporting our theoretical findings.

\textbf{Notation.} We conclude this section by introducing some notational conventions that will be used throughout the paper. The all-$1$ and all-$0$ column vectors are denoted as $\vect{1}$ and $\vect{0}$, respectively, with appropriate dimension depending on the context.  We denote the identity matrix by ${I}$, where the dimension will be clear from the context. 
	A matrix $W$ is row stochastic when its entries are non-negative and $W\vect{1} = \vect{1}$. A matrix $W$ is said to be Schur stable if its spectral radius, which coincides with the absolute value of its maximum eigenvalue, is less then 1, i.e. $\rho(W)=|\lambda_{max}(W)|<1$.  Given a matrix $X\in \mathbb{R}^{\mc V\times\mc V}$ where $\mc V$ is a finite set and given $\mc A, \mc B\subseteq \mc V$, we indicate with $X_{\mc{AB}}$ the submatrix of $X$ having rows in $\mc A$ and columns in $\mc B$. We use the notation $X_{\mc B}$ for $X_{\mc V \mc B}$. Moreover, given a subset $\mc A$ of $\mc V$, we define $X_{- \mc A} = X_{\mc V \backslash \mc A}$. The notation $A = [\beta]$, with $\beta \in \mathbb{R}^n$ will refer to a $n\times n$ diagonal matrix with $A_{ii}=\beta_i$. We will use the notation $C \succ 0$ to indicate a positive definite matrix. 

	\section{Model and Problem Formulation}

We model the social network as a finite undirected weighted graph  $\mc G = (\mc V, \mc E, W)$, where the node set $\mc V$  represents agents and the link set $\mc E$ represents interactions. Social influence among the agents is measured by the link weights $W_{ij}=W_{ji}$ that are collected in a symmetric weight matrix $W=W'$ in $\mathbb{R}_{+}^{\mc V \times \mc V}$, {which is assumed to be known}. Social influence weights are such that $W_{ij}>0$ if and only if $\{i,j\} \in \mc E$. 

We consider a DeGroot opinion dynamics model with stubborn agents and noise, as illustrated below. 
Every agent $i$ in $\mc V$ is endowed with a time-varying state $X_i(t)$ representing his opinion  on a particular topic of interest. 
Let $\mc S \subset \mc V$ be a globally reachable subset of agents, to be referred as stubborn agents. We shall refer to the set $\mc R = \mc V \backslash \mc S$ of the non-stubborn agents as regular agents. 
Let $P$ in $\mathbb{R}^{\mc V \times \mc V}$  be  the row-stochastic normalized weight matrix with entries $P_{ij}= \frac{W_{ij}}{w_i}$ where $w_i=\sum_{j \in \mc V} W_{ij}$. 	We then assume that individuals update their opinion by taking into account the views of others according to the following update rule:
\be \label{eqn:DeGroot general}
\begin{array}{lll}
	X_i(t+1) &= \sum_{j \in \mc V}P_{ij} X_j(t)+V_i(t+1) \qquad &\text{if } i \in \mc R \\
	X_i(t+1) &=u_i \qquad &\text{if } i \in \mc S
\end{array}	
\ee
where, for every stubborn agent $i$ in $\mc S$, $u_i $ is his constant opinion, whereas, for every regular agent $j$ in $\mc R$,  $V_i(t)$ is a sequence of i.i.d.~random noise variables with zero mean and bounded second moment $\sigma_i^2$. Throughout, we shall assume that the noise variables $V_{i}(t)$ are independent across the regular agents $i$ in $\mc R$. Let $V(t)=(V_i(t))_{i\in\mc R}$ be the noise vector, then 
$$\mathbb E[V(t)]=0\,,\qquad \Sigma=\mathbb{E}[V(t)V(t)']=[\sigma^2]\,.$$
{The regular agents gets influenced through the network $\mc G$ and, at every time instant, update their opinion to a noisy average of their neighbors'. Conversely, the stubborn one do not interact and maintain constant their opinion.}
{We now focus on the opinion of the regular agents and, to this aim, we gather all their opinions at time $t$ in a vector $X(t):=(X_i(t))_{i \in \mc R}$.}
We can thus rewrite the dynamics of the regular agents' opinions as:
\be \label{eqn:dynamical system}
X(t+1)=AX(t)+Bu+V(t+1)
\ee
where $A:=P_{\mc R \mc R}$ refers to the influence of regular agents; $u:=(u_i)_{i \in \mc S}$ is the vector of stubborn nodes opinions, which is assumed {to be} deterministic; $B:=P_{\mc R \mc S}$ indicates their impact on regular agents.
In particular, it follows from Theorem 1 in \cite{Ravazzi2015_noiseconvergence}   that $X(t)$ converges in distribution to a random equilibrium opinion vector $X$. 
\footnote{
	Recall that a sequence of random vectors $X(t)$ taking values in $\R^n$ converge in distribution to a random vector $X$ if $\mathbb E[\varphi(X(t))]\stackrel{t\to+\infty}{\longrightarrow}\mathbb E[\varphi(X)]$ for every continuous bounded function $\varphi:\R^n\to\R$. 
}

Now, assume we are interested in estimating the mean equilibrium opinion {of regular agents, i.e.}
$$Y = \frac{1}{n}\sum_{i \in \mc R}X_i,$$
but that we can observe only a limited number of agents' opinions. Specifically, we aim to select a subset $\mc K \subseteq \mc R$ of at most $s$ regular  agents and a function $g:\mathbb{R}^{\mc K} \rightarrow \mathbb{R}$ in order to minimize the mean squared error $$ \mathbb{E}[(Y-g(X_{\mc K}))^2 ]\,.$$
Given $\mc K \subseteq \mc R$, it is well known (see Section 1.5.1 in \cite{ross1996stochastic}) that the best estimator $g^*(X_{\mc K})$ of $Y$ that minimizes the mean squared error is  
$$g^*(X_{\mc K}) = \mathbb{E}[Y|X_{\mc K}]\,,$$ 
i.e., 
$$G(\mc K):= \min_{g(X_{\mc K})} \mathbb{E}[(Y-g(X_{\mc K}))^2 ]=\mathbb{E}[(Y-g^*(X_{\mc K}))^2]\,. $$ 
Using this result we retrieve that
\begin{equation}\label{MSE}
	\begin{array}{rcl}
		G(\mc K) &=& \mathbb{E}[(Y-g^*(X_{\mc K}))^2 ]\\
		&=&	\mathbb{E}[\mathbb{E}[(Y-\mathbb{E}[Y|X_{\mc K}])^2|X_{\mc K}]] \\
		&=& \mathbb{E}[\Var(Y|X_{\mc K})]	
	\end{array}	
\end{equation}
i.e., the mean squared error coincides with the residual variance on $Y$ given observation $X_{\mc K}$. 

Our goal is to select the regular agents subset  $\mc K \subseteq \mc R$ of given cardinality $s$ to observe in order to best estimate $Y$, i.e.,   
\be \label{problem:minMSE} 	{\mc K}^* = \underset{\ds |\mathcal{K}|=s}{\arg \min} \; G(\mathcal{K}) \,.\ee 
Let us observe that in many situations it may be convenient to work with variance reduction as objective function, i.e.
\be \label{eqn:MSE-VarRed} F(\mc K) := \Var(Y)-G(\mc K) \,. \ee
This leads us to an alternative formulation of the previously introduced problem which consists in the maximization of $ F(\mc K)$, i.e.
\be	\label{problem:maxF}
{\mc K^*} = \underset{\ds |\mathcal{K}|=s}{\arg \max} \; F(\mathcal{K}) \,.
\ee

\section{Preliminary Results}
Let us now present some preliminary results concerning the equilibrium opinion vector $X$ which will be usefull in the following.  All the detailed proofs are reported in Appendix. 
\begin{proposition} \label{prop:general covariance}
	Let $X$ be the random equilibrium opinion vector of the DeGroot opinion dynamics  \eqref{eqn:dynamical system} on an undirected graph  with nonempty globally reachable set of stubborn nodes. Then, $X$ has expected value
	\be \label{eq:mean_submodularCase} \mathbb{E}[X] =(I-A)^{-1}Bu,,\ee
	and covariance matrix	\be \label{eq:cov_submodularCase} 	C=\mathbb{E}[(X-\mathbb{E}[X])(X-\mathbb{E}[X])']=\Sigma(I-A^2)^{-1}
	\ee
\end{proposition}
\begin{proof}
	See Appendix-\ref{app:Proposition1}.
\end{proof}
\begin{corollary} \label{lemma: C properties}
	The covariance matrix $C:=\mathbb{E}[(X-\mathbb{E}[X])(X-\mathbb{E}[X])']$ satisfies the following properties:
	\begin{itemize}
		\item[(i)] C is positive definite, i.e. $C \succ 0$;
		\item[(ii)] the precision matrix $H:= C^{-1}$ is symmetric positive definite and such that  $H_{ij} \leq 0$ for any $i \neq j$.		
	\end{itemize} 
\end{corollary}
\begin{proof}
	See Appendix-\ref{app:Proposition1}.
\end{proof}
%
We can retrieve now a general formulation  for the objective functions of interest. This adapts previous known results from the literature \cite{KempeApproxSubmod} to our specific setting. For completeness a simple detailed proof is shown in Appendix.
\begin{proposition}
	\label{prop:F formulation}
	Given arbitrary $\mc K \subseteq \mc R$, 
	\begin{itemize}
		\item [(i)] the variance reduction is computed as	\begin{equation}
			\label{function_general}
			F(\mc K) =  (C \vect 1)'_\mc K(C_{\mc K \mc K})^{-1}( C {\vect 1})_\mc K \,,
		\end{equation}
		where $C=\Sigma(I-A^2)^{-1}$.
		\item [(ii)] the mean squared error is computed as \be \label{eqn:MSPE formulation}
		G(\mc K) = \vect 1_{-\mc K}'(H_{-\mc K-\mc K})^{-1}\vect 1_{-\mc K}\,.
		\ee
		where $H = (I-A^2)\Sigma^{-1}.$
	\end{itemize}
\end{proposition}
\begin{proof}
	See Appendix-\ref{app:Proposition 2}.
\end{proof}

As a direct consequence of Proposition \ref{prop:F formulation}, when $|\mc K|=1$ we can introduce a simplified expression for $F(\{k\})$.

\begin{corollary} \label{cor:F single obs}
	If $\mc K =\{k\}$, variance reduction $F(\{k\})$ is computed as
	\[
	F(\{k\}) = \sigma_k \eta_k,
	\]
	where $\sigma_k = \sqrt{\Sigma_{kk}}$ and  $\eta_k:= \frac{((I-A^2)^{-1}\vect{1})_k^2}{((I-A^2)^{-1})_{kk}}$\,.
\end{corollary}
Analysing the function in Corollary \ref{cor:F single obs} we can highlight an interesting connection with a known centrality measure, the intercentrality, introduced by \cite{keyPlayer} for the well-known Key Player problem. 
\begin{remark}
	Let us preliminary define a new graph $\hat{\mc G} = (\hat{\mc V}, \hat{\mc E})$	s.t.  $\hat{\mc V} = \mc V$ and $\hat{\mc E}$ defined s.t. $(i,j) \in \hat{\mc E}$ iff it exists a path of lenght 2 connecting $i$ and $j$ in $\mc G$. $\hat {\mc G}$ is known as the the 2-hop neighbour graph \cite{Jin2006_multihop} associated to $\mc G$. 
	We notice that the function $\eta_k$, introduced in Corollary \ref{cor:F single obs}, coincides with the intercentrality \cite{keyPlayer} of agent $k$, defined on the 2-hop neighbour graph $\hat{\mc G}$.
	
	For the sake of clarity, in Fig.\ref{fig:2hopcycle} the 2-hop neighbour graph associated to a cycle graph with 7 nodes is represented. More precisely, the links in $\hat{\mc E}$ are the one highlighted in red in Figure and connect those nodes which are at distance 2 in the original graph $\mc G$. 
	\begin{figure}[!h]
		\centering
		\includegraphics[width=0.6\linewidth]{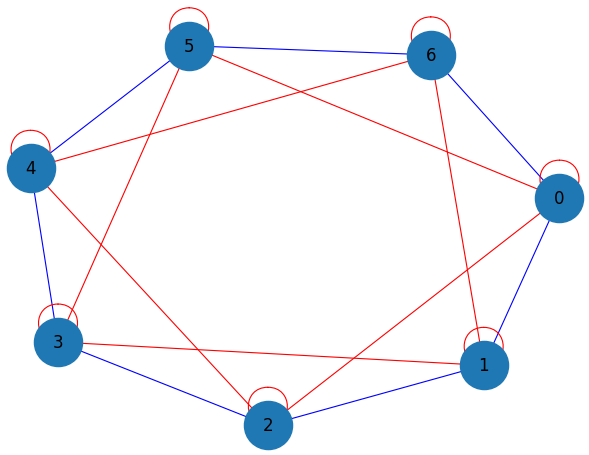}
		\caption{Cycle graph with 7 nodes. The blue links represent $\mc E$, while the red one are associated to $\hat{\mc E}$.}
		\label{fig:2hopcycle}
	\end{figure}
	
\end{remark}
Finally, referring to our variance reduction function $F(\{k\})$, we highlights also the distinction between our concept of the ’most informative’ node and what is typically referred to in the literature as the ’most influential’ node, which is commonly referred to as the node maximizing the Bonacich centrality. 
{This difference arises because Bonacich centrality is based on the total number of walks involving each node, whereas in variance reduction, this term affects only the numerator of the ratio introduced in Corollary \ref{cor:F single obs}. The denominator, which accounts for the number of cycles, properly rescales the expression, leading to the observed distinction. }
Analogous result has been shown in our previous work \cite{IFAC} in the particular case of Gaussian distributed stubborn nodes and in absence of additional noise. The following example corroborates this key difference also in this different setting of interest.
\begin{example}
	Consider the Watts-Strogatz graph with $|\mc V|=15$ and 3 stubborn nodes, represented in Fig.\ref{fig:WS}.  We compare the results of the proposed Variance Reduction measure $F(\{k\})$ with the Bonacich centrality. The colormap reflects the centrality measure for each node in the network (notice that in Figure the values are normalized). It is worth noting that our method diverges from existing approaches: the node identified as the most informative (green node in Fig.\ref{fig:WS}a) does not coincide with the one with highest Bonacich centrality (green node in Fig.\ref{fig:WS}b ).
	\begin{figure}[!h]
		\centering
		\subfloat[Variance Reduction]{
			\includegraphics[width=0.48\linewidth]{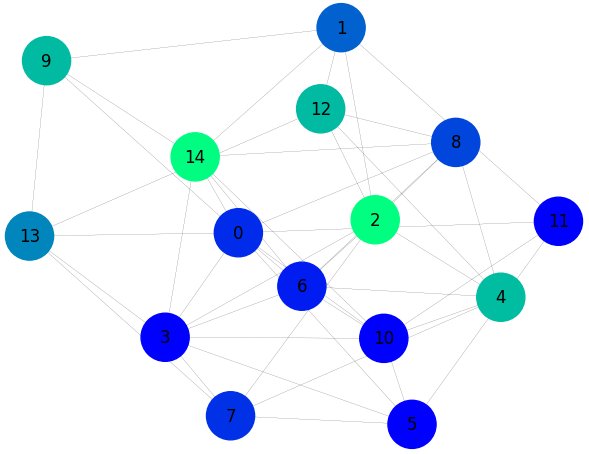}   \label{fig:WS_VarRed}} 
		\subfloat[Bonacich]{
			\includegraphics[width=0.48\linewidth]{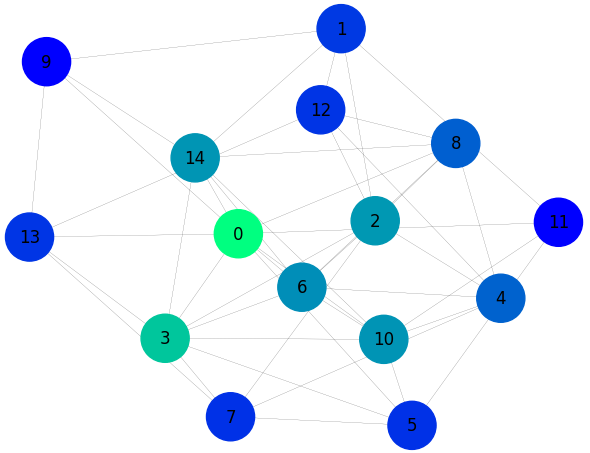}   \label{fig:WS_Bonacich}}
		\caption{Comparison of centrality measures over a Watts-Strogatz graph. The colormap indicates the normalized Variance Reduction and the normalized Bonacich Centrality, respectively..}
		\label{fig:WS}
	\end{figure}
	
\end{example}
\section{Main Results }
In this section, we present our main result, stating that, {in the considered setting, with no extra assumption on the model parameters,} the variance reduction function $F(\mc K)$ is sub-modular, equivalently, that the minimal mean-square error function $G(\mc K)$ is super-modular, as per the following definition. 



\begin{definition}
	Given a set $\mc Z$, a function $f:2^{\mc Z} \rightarrow \mathbb{R}$ is 
	\begin{enumerate}
		\item[(i)]
		sub-modular if 
		\be\label{eqn:def_submodularity}
		f(\mathcal{A}\cup \{k\}) - f(\mathcal{A}) \geq f(\mathcal{B}\cup \{k\})-f(\mathcal{B})\,,\end{equation}
	for every $\mathcal{A}\subseteq\mathcal{B}\subseteq\mathcal{V}$ and $k \in \mathcal{V}\backslash\mathcal{B}$;
	\item[(ii)]
	super-modular if $-f$ is sub-modular, i.e., if 
	\be\label{eqn:def_supermodularity}
	f(\mathcal{A}\cup \{k\}) - f(\mathcal{A}) \leq f(\mathcal{B}\cup \{k\})-f(\mathcal{B})\,,\end{equation}
for every $\mathcal{A}\subseteq\mathcal{B}\subseteq\mathcal{Z}$ and $k \in \mathcal{Z}\backslash\mathcal{B}$.
\end{enumerate}
\end{definition}

\begin{theorem} \label{th:submod}
Let $X$ be the random equilibrium opinion vector of the DeGroot opinion dynamics on an undirected graph  with nonempty globally reachable set of stubborn nodes. 
For every subset of regular nodes $\mc K\subseteq\mc R$, let $G(\mc K)$ and $F(\mc K)$ be, respectively, the minimal mean-square error \eqref{eqn:MSPE formulation} and the maximal variance reduction \eqref{function_general} in the estimation of the mean opinion $Y$ from $X_{\mc K}$.
Then, 
\begin{enumerate}
\item[(i)] $F:2^{\mc R}\to\R$ is {sub}-modular; 
\item[(ii)] $G:2^{\mc R}\to\R$ is {super}-modular. 
\end{enumerate}
\end{theorem}
\begin{proof} 
We will just prove point (i) since that, together with \eqref{eqn:MSE-VarRed}, readily implies point (ii). 
Let $D$ in $\R^{\mc R\times\mc R}$ have entries  $$D_{ii}=1/\sqrt{H_{ii}}\,,\qquad D_{ij}=0\,,\qquad \forall i\ne j\in\mc R\,.$$ 
Let $H=(I-A)^2{\Sigma^{-1}}$ be the precision matrix of $X$  and let
$$\bar{H} = DH D\,,\qquad 
\Theta = I-\bar{H}\,.$$ 
Now, observe that, on the one hand 
$$\Theta_{ii}=1-\overline H_{ii}=1-D_{ii}H_{ii}D_{ii}=1-\frac{H_{ii}}{\sqrt{H_{ii}H_{ii}}}=0\,,$$ for every $i$ in $\mc R\setminus\mc K$, on the other hand 
$$\Theta_{ij}=-\overline H_{ij}=-D_{ii}H_{ij}D_{ij}=\frac{-H_{ij}}{\sqrt{H_{ii}H_{jj}}}\ge0\,,$$
for every $i\ne j$ in $\mc R$, where the inequality follows from Corollary \ref{lemma: C properties}(i).  Hence, $\Theta$ is a symmetric matrix with nonnegative entries.  
The Perron-Frobenius Theorem then implies that $\Theta_{-\mc K,-\mc K}$ has a real dominant eigenvalue $\lambda_{\mc K}\ge0$ such that $\rho(\Theta_{-\mc K-\mc K})=\lambda_{\mc K}$. Moreover, since $H$ is symmetric positive definite by Corollary \ref{lemma: C properties}(i), so is $\bar{H} = DH D$. By the Sylvester criterium, also 
$\bar{H}_{-\mc K -\mc K}=I-\Theta_{-\mc K-\mc K}$ is symmetric positive definite, so that necessarily $\lambda_{\mc K}<1$. 
It follows that $\rho(\Theta_{-\mc K-\mc K})<1$, so that  $\bar{H}_{-\mc K -\mc K}$ is invertible and that its inverse coincides with the limit of the geometric series
\be\label{barH}(\bar{H}_{-\mc K-\mc K})^{-1} =(I-\Theta_{-\mc K-\mc K})^{-1} =  \sum_{l=0}^{+\infty} (\Theta_{-\mc K-\mc K})^l\,.\ee 
For every $i$ and $j$ in $\mc R\setminus\mc K$, and $l\ge0$, let 
$$\Gamma^{l}_{ij}(\mc K)=\left\{\!\gamma\!=\!(\gamma_0,\!\gamma_1,\!\ldots,\!\gamma_l)\in(\mc R\setminus\mc K)^{l+1}\!\!:\gamma_0=i,\gamma_l=j\right\}\,,$$
and define
$$\theta_{\gamma}=\prod_{1\le h\le l}\Theta_{\gamma_{h-1}\gamma_h}\,,\qquad\forall\gamma\in\Gamma^{l}_{ij}(\mc K)\,,$$
(with the convention that $\theta_{\gamma}=1$ when $l=0$ and $\gamma=(i)$). 
Then, \eqref{barH} reads 
\be\label{barHij}
((\bar{H}_{-\mc K-\mc K})^{-1})_{ij} = \sum_{l=0}^{+\infty} \sum_{\gamma\in\Gamma^{l}_{ij}(\mc K)}\theta_{\gamma}
\,.\ee
Now, consider an arbitrary $k$ in $\mc R\setminus\mc K$ and let $\ov{\mc K}=\mc K\cup\{k\}$. Then, for every $i$ and $j$ in $\mc R\setminus\ov{\mc K}$,
\be\label{barHij-}
((\bar{H}_{-\mc K-\mc K})^{-1})_{ij} -((\bar{H}_{-\ov{\mc K}-\ov{\mc K}})^{-1})_{ij} 
=\sum_{s=0}^{+\infty} \sum_{\substack{\gamma\in\Gamma^{l}_{ij}(\mc K)\\k\in\{\gamma_1\,\ldots,\gamma_{n-1}\}}}\theta_{\gamma}\,.\ee

Now, observe that  \eqref{eqn:MSPE formulation} implies that 
\be \label{eqn:G_form}
\ba{rcl}	
G(\mc K) &=&  \vect 1_{-\mc K}'(H_{-\mc K-\mc K})^{-1}D\vect 1_{-\mc K}\\
&=& \vect 1_{-\mc K}'((D^{-1}\ov HD^{-1})_{-\mc K-\mc K})^{-1}\vect 1_{-\mc K}\\
&=& \vect 1_{-\mc K}'(D^{-1}_{-\mc K-\mc K}\ov H_{-\mc K-\mc K}D^{-1}_{-\mc K-\mc K})^{-1}\vect 1_{-\mc K}\\
&=& \vect 1_{-\mc K}'D_{-\mc K-\mc K}(\ov H_{-\mc K-\mc K})^{-1}D_{-\mc K-\mc K}\vect 1_{-\mc K}\\
&=& 
\ds\sum_{i,j\in\mc R\setminus\mc K}\frac{((\bar{H}_{-\mc K-\mc K})^{-1})_{ij}}{\sqrt{H_{ii}H_{jj}}}\,.
\ea\ee
It then follows from \eqref{eqn:MSE-VarRed} and \eqref{eqn:G_form} that 
$$
\ba{rcl}F(\ov{\mc K})\!-\! F(\mc K) \!\!\!\! \!\!\!
&=&G(\mc K)-G(\ov{\mc K})\\
&=&\ds\sum_{i,j\in\mc R\setminus\mc K}\frac{((\bar{H}_{-\mc K-\mc K})^{-1})_{ij}}{\sqrt{H_{ii}H_{jj}}}\\
&&-\ds\sum_{i,j\in\mc R\setminus\ov{\mc K}}\frac{((\bar{H}_{-\ov{\mc K}-\ov{\mc K}})^{-1})_{ij}}{\sqrt{H_{ii}H_{jj}}}\,.\\
&=&E_1(\mc K)+E_2(\mc K)+E_3(\mc K)\,,\ea
$$
where, using  \eqref{barHij}, and \eqref{barHij-} we have that 
$$E_1(\mc K)=\frac{((\bar{H}_{-\mc K-\mc K})^{-1})_{kk}}{H_{kk}}=\sum_{l=0}^{+\infty} \sum_{\substack{\gamma\in\Gamma^{l}_{kk}(\mc K)}}\frac{\theta_{\gamma}}{H_{kk}}$$
$$\ba{rcl}E_2(\mc K)
&=&\ds2\sum_{i\in\mc R\setminus\ov{\mc K}}\frac{((\bar{H}_{-\mc K-\mc K})^{-1})_{ki}}{\sqrt{H_{ii}H_{kk}}}\\
&=&\ds2\sum_{i\in\mc R\setminus\ov{\mc K}}\sum_{l=0}^{+\infty} \sum_{\substack{\gamma\in\Gamma^{l}_{ik}(\mc K)}}\frac{\theta_{\gamma}}{\sqrt{H_{ii}H_{kk}}}\,,\ea$$
$$\ba{rcl}E_3(\mc K)&=&\ds\sum_{i,j\in\mc R\setminus\ov{\mc K}}\frac{((\bar{H}_{-\mc K-\mc K})^{-1})_{ij}}{\sqrt{H_{ii}H_{jj}}}\\
&=&\ds\sum_{i,j\in\mc R\setminus\ov{\mc K}}\sum_{l=0}^{+\infty} \sum_{\substack{\gamma\in\Gamma^{l}_{ik}(\mc K)}}\frac{\theta_{\gamma}}{\sqrt{H_{ii}H_{jj}}}\,.\ea$$
Notice that, since $\theta_{\gamma}\ge0$, for every $\mc A\subseteq\mc B\subseteq\mc R$, we have that 
$$E_h(\mc A)\ge E_h(\mc B)\,,\qquad \forall h=1,2,3\,,$$  
so that 
$$\ba{rcl}F(\mc A\cup\{k\})-F(\mc A)
&=&E_1(\mc A)+E_2(\mc A)+E_3(\mc A)\\ 
&\ge&E_1(\mc B)+E_2(\mc B)+E_3(\mc B)\\ 
&=&F(\mc B\cup\{k\})-F(\mc B)\,,\ea
$$
thus showing that $F$ is a sub-modular function. 
\end{proof}

\section{Greedy Algorithm Design}
The result of Theorem \ref{th:submod} enhances the applicability of a simplified approach to solve the combinatorial optimization probem in \eqref{problem:maxF}. As well-established in the literature  \cite{Nemhauser}, let us in particular define a Greedy Algorithm adapted to our problem:
\begin{enumerate}
\item Denote with $\mc K_t$ the subset of nodes selected at $t$-th algorithm iteration and initialize $\mc K_0=\emptyset$. 
\item For each iteration $t=1, \dots, s$, with $s$ cardinality of the observed set to select, set $\mc K_t = \mc K_{t-1} \cup \{i_t\}$ with
$$ i_t \in \arg\max_{i \in \mc R \backslash \mc K_{t-1}} F(\mc K_{t-1}\cup \{i_t\}) \,.$$
\item At the end of $s$-th iteration, the algorithm will propose $\mc K_s$ as an approximated solution to problem \eqref{problem:maxF}. 
\end{enumerate}
Despite the sub-optimality of the proposed solution, from \cite{Nemhauser} we have some guarantees on method performances, as properly stated in the following Corollary.
\begin{corollary} \label{cor:greedy_opt}
The Greedy Algorithm is a $(1-1/e)$-approximation for the problem of maximizing the variance reduction introduced in \eqref{problem:maxF}, more precisely, fixed the cardinality of the observable set to $s$ and given $\mc K_s$ the solution proposed by the algorithm, it holds
$$ F(\mc K_s) \geq \left(1-\frac{1}{e}\right) \max_{\mc K \subseteq \mc R: |\mc K|=s}F(\mc K)$$
\end{corollary}

{We notice the complexity reduction intrinsic in the greedy algorithm where, at every step, we need to solve a maximization problem on a set whose cardinality is bounded by $n$. In this way, only $O(ns)$ evaluations of the variance reduction function are needed, instead of the $O(n^s)$ needed in the original problem.}


The computational cost can be further reduced by optimizing the evaluation of the $F$ function, specifically by avoiding matrix inversion at each iteration, which typically requires $O(s^3)$ operations.  In the following, we focus on deriving an alternative formulation of $F$, which will enable the development of a more efficient greedy algorithm. The results presented follows from the application of linear algebra tools related to matrices such as Schur Complement and Shermann Morrison formula. For detailed proofs see Appendix-\ref{app:Greedy}.
\begin{proposition} \label{prop:iterative step}
	Let $\mc K$ be the observed set at step $t-1$ and $i$ be the candidate extra node to observe at $t$-iteration. Given the variance reduction function $F$, the iterative step is defined as 
	\[
	F(\mc K \cup \{i\}) = \delta + \frac{(C\alpha)_i^2}{C_{ii}}\left[1+C_{i \mc K}(C_{\mc K \mc K})^{-1}\left(1-2\frac{(C\alpha)_{\mc K}}{C_{ii}}\right)\right]
	\]
	where $\delta=(C\alpha)_{\mc K}'(C_{\mc K \mc K})^{-1}(C\alpha)_{\mc K}$ is a constant which depends only on $\mc K$.
\end{proposition}
\begin{proposition} \label{cor:iterative_inverse}
	The inverse matrix $(C_{\mc K \mc K})^{-1}$ can be computed iteratively without requiring any matrix inversion. More precisely, 
	\begin{itemize}
		\item for $t=1$ since $\mc K_1 = \{i\}$ then $(C_{\mc K_1 \mc K_1})^{-1}=\frac{1}{C_{ii}}$
		\item 	 In the subsequent steps, indicated with $I$ the inverse matrix at time $t-1$, i.e. $(C_{\mc K_{t-1} \mc K_{t-1}})^{-1}$, then $I = (C_{\mc K_{t} \mc K_{t}})^{-1}$ is computed as
		\[
		I = \begin{bmatrix}
			\bar I & -\bar I \frac{C_{\mc K i^*}}{C_{i^* i^*}}\\
			-\frac{C_{\mc K i^*}}{C_{i^* i^*}}\bar I & \frac{1}{C_{i^* i^*}}+ \frac{C_{i^* \mc K}}{C_{i^* i^*}}IC_{\mc K i^*}
		\end{bmatrix}
		\]
		where $\bar I =  I+\frac{IC_{\mc K i^*}C_{i^* \mc K}I}{C_{i^*i^*-C_{i^*\mc K}IC_{\mc K i^*}}}$
	\end{itemize}
	
\end{proposition}

Based on Proposition \ref{prop:iterative step} - \ref{cor:iterative_inverse}, we formally define the pseudo-code for the greedy algorithm, as presented in Algorithm  \ref{alg:greedy}.

\begin{algorithm} [!h]
	\caption{Greedy algorithm for optimal subset selection}
	\KwData{$C$, $s$, $\mc R$   }
	\KwResult{$\mc K$  }
	$\mc K \gets \emptyset$\;
	$F_{max} \gets 0$\;
	$I \gets 0$\;
	\For{$t=1, \dots, s$}{
		$F_{max} \gets 0$\;
		\For{$i \in \mc R \backslash \mc K$}{
			$F_i \gets \frac{(C\alpha)_i^2}{C_{ii}}\left[1+C_{i \mc K}I \left(1-2\frac{(C\alpha)_{\mc K}}{C_{ii}}\right)\right] $\;
			\If{$F_i > F_{max}$}{
				$F_{max}\gets F_i$\;
				$i^*\gets i$\;
			}
		}
		$\mc K \gets \mc K \cup \{i^*\}$\;
		$\bar I \gets I+\frac{IC_{\mc K i^*}C_{i^* \mc K}I}{C_{i^*i^*-C_{i^*\mc K}IC_{\mc K i^*}}}$\;
		$I \gets \begin{bmatrix}
			\bar I & -\bar I \frac{C_{\mc K i^*}}{C_{i^* i^*}}\\
			-\frac{C_{\mc K i^*}}{C_{i^* i^*}}\bar I & \frac{1}{C_{i^* i^*}}+ \frac{C_{i^* \mc K}}{C_{i^* i^*}}IC_{\mc K i^*}
		\end{bmatrix}$\;
	}
	\label{alg:greedy}
\end{algorithm}
We now analyze the computational complexity of the proposed greedy algorithm in relation to the exact brute-force approach presented in the previous section, to remark its scalability also for large-scale networks.
\begin{proposition}[Computational Cost] \label{prop:comp_cost_greedy}
	The greedy approximation proposed in Algorithm \ref{alg:greedy} achieves an overall computational cost of $O(ns^3)$ for the subset selection problem, significantly reducing the complexity compared to the exact brute-force algorithm which requires $O(n^s s^3)$ operations.
\end{proposition}
\begin{proof} Let us analyse the computational cost of the two main steps:
	\begin{itemize}
		\item[(i)] Since one element is added at each step the total number of iterations of the algorithm is equal to $s$ and the total number of evaluation of the objective function $F(\mc K)$ are thus equal to $O(ns)$;
		\item[(ii)] Applying the iterative computations introduced in Proposition \ref{cor:iterative_inverse} the computational cost of matrix inversion is reduced to $O(s^2)$ and this coincides with the cost associated to the single evaluation of the objective function $F(\mc K)$.
	\end{itemize}
On the other side, focusing on the exact brute-force algorithm, it requires a combinatorial complexity for the subset selection, i.e. $O(n^s)$, combined with the cost of matrix inversion, i.e. $O(s^3)$.
	Collecting these results the thesis follows. 
	
\end{proof}
\color{black}
Let us now show two simple simulations  to support the validity of this method. 

\begin{example}[Watts-Strogatz Network]
Let us consider a randomly generated Watts-Strogatz network with 15 nodes, 3 of which are stubborn. Suppose that we are interested in observing up to 4 nodes and compare the results obtained both using the exact computation from equation \eqref{function_general} and applying the Greedy Algorithm. In Fig.\ref{fig:wsgreedy15n} is represented the percentage of residual variance conditioned to the observation of subset $\mc K$. We highlight how the Greedy Algorithm achieves comparable performances in terms of residual variance while significantly reducing the computational cost. To be thorough, in Fig.\ref{fig:WS_selection} it is reported the proposed selection for a 4 nodes observation set both in case of exact computation and Greedy Algorithm application.
\begin{figure}[!h]
\centering
\includegraphics[width=0.7\linewidth]{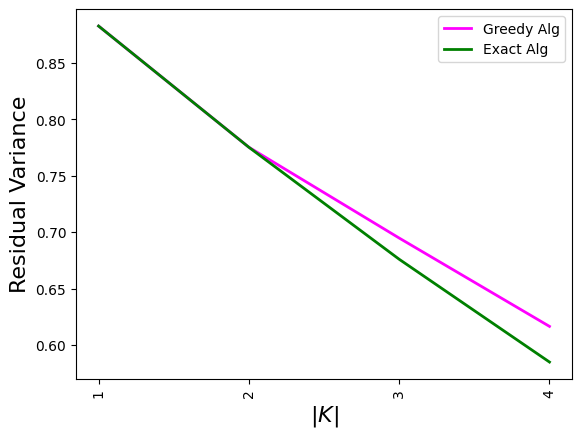}
\caption{Comparison between the percentage of residual variance gained through the exact method vs the Greedy Algorithm, as function of the cardinality of the observed set $|\mc K|$, for a Watts-Strogatz Network.}
\label{fig:wsgreedy15n}
\end{figure}
\begin{figure}[!h]
\centering
\subfloat[Exact Selection]{
	\includegraphics[width=0.48\linewidth]{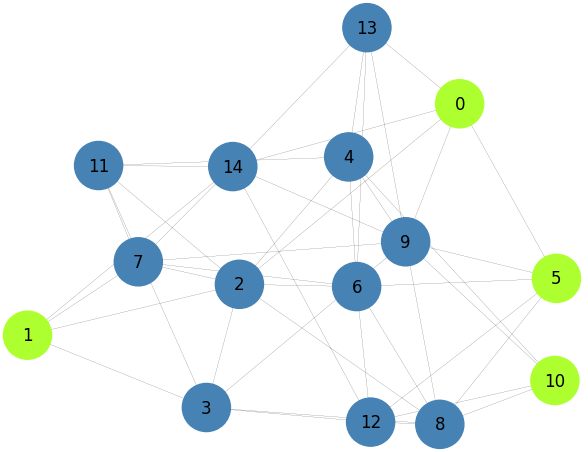}   \label{fig:WS_Exact}} 
\subfloat[Greedy Selection]{
	\includegraphics[width=0.48\linewidth]{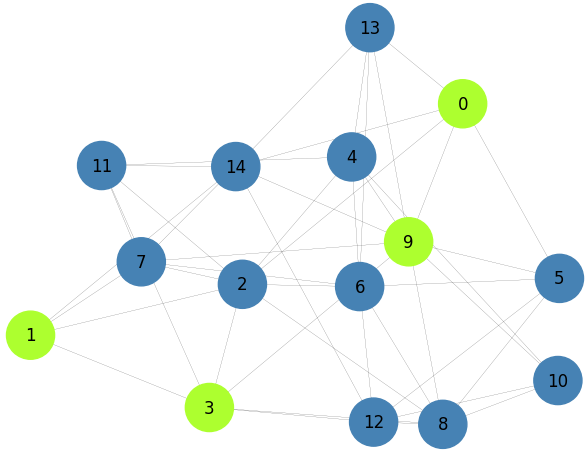}   \label{fig:WS_Greedy}}
\caption{Comparison between the nodes selection through  exact method vs Greedy Algorithm, fixed cardinality $s=4$, for a Watts-Strogatz Network.}
\label{fig:WS_selection}
\end{figure}
\end{example}
\begin{example}[Real-case scenario]
\begin{figure}[!h]
\centering
\includegraphics[width=0.7\linewidth]{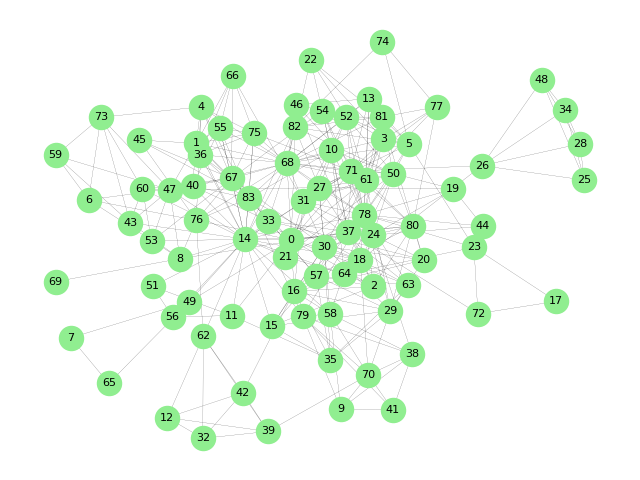}  
\caption{Malawi network used in the case study
	analized.  }
\label{fig:Net_RealDB}
\end{figure}
Here we demonstrate our approach on a real-case scenario.
In particular, we apply our methodology to the network associated to social contacts in a village of rural Malawi \cite{malawi_net}, whose
dataset is available in Sociopatterns. The network obtained has 86 individuals, {347} weighted undirected links, and it is
illustrated in Fig. \ref{fig:Net_RealDB}.
Despite the worse performances, highlighted in Fig. \ref{fig:Greedy_realDB}, achieved by the greedy algorithm, its computational cost makes problems of larger dimension approachable.	

\begin{figure}[!h]
\centering
\includegraphics[width=0.7\linewidth]{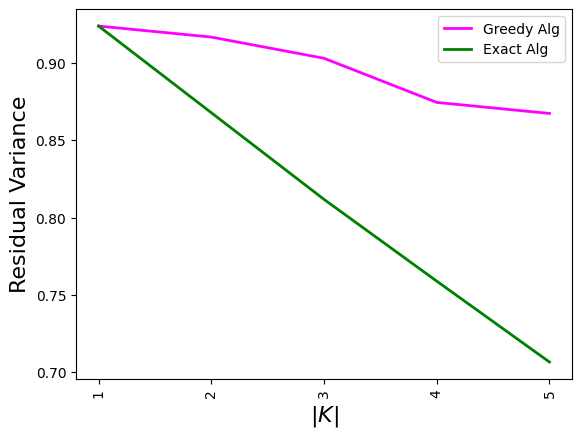}
\caption{Comparison between the percentage of residual variance gained through the exact method vs the Greedy Algorithm, as function of the cardinality of the observed set $|\mc K|$, for the Malawi Network.}
\label{fig:Greedy_realDB}
\end{figure}
\end{example}
\section{Conclusions}
In this work we have analysed a subset selection problem for a particular Opinion Dynamics model with deterministic stubborn agents and random external noise. First, we successfully extended our previous results in \cite{IFAC}, which focused solely on the case of an initial Gaussian distribution and in absence of external noise, to this more general and relevant setting. We thus showed that the formulation of the variance reduction function remains unchanged also in case of arbitrary distribution. Additionally, we proved that the variance reduction function $F(\mc K)$, objective function in the optimization problem, exhibits the submodularity property in this specific setting. Consequently, based on well-known results in the literature, we defined an appropriate applicable greedy algorithm which identifies a good sub-optimal solution to the problem. Some simple simulations have been presented to support the theoretical findings.

The results obtained, along with those previously presented in \cite{IFAC}, hence lay the groundwork for a more general and comprehensive problem formulation and analysis aimed at addressing both the stochasticity of stubborn agents and the occurrence of additional random noise in transmission, which will be further explored in future work. {Moreover, a challenging problem to address could be to evaluate how noisy or partial information on network structure influence the objective function of interest. }

\appendix
\subsection{Proof of Proposition \ref{prop:general covariance} and Corollary \ref{lemma: C properties} }\label{app:Proposition1}
\begin{proof}[Proof of Proposition \ref{prop:general covariance}]
	First notice that $A$ is a Schur stable matrix since it is  the restriction of the stochastic matrix $P$ to the complement of a globally reachable node subset  $\mc S\subseteq\mc V$ (c.f.~\cite{Como.Fagnani:2016}). Hence, $I-A$ is invertible. 
	
	Since $X(t)$ converges in distribution to $X$, then $$X\sim AX+Bu+V\,,$$ where $\sim$  stands for equality in distribution and $V\sim V(t)$ is an independent copy of the noise vector.  
	Since $\mathbb{E}[V]=0$, we have that 
	$$ \mathbb{E}[X]=A\mathbb{E}[X]+Bu\,.$$
	The above implies that $(I-A)\mathbb{E}[X]=Bu$, from which \eqref{eq:mean_submodularCase} follows, thanks to the invertibility of $I-A$.

	On the other hand 
	$$\ba{rcl}C
	&=&\ds\mathbb{E}[(X-\mathbb{E}[X])(X-\mathbb{E}[X])']\\
	&=&\ds\mathbb{E}[(A(X-\mathbb{E}[X])+V)(A(X-\mathbb{E}[X])+V)']\\	
	&=&\ds A\mathbb{E}[(X-\mathbb{E}[X])(X-\mathbb{E}[X])']A'+\mathbb E[VV']\\
	&=&\ds ACA'+\Sigma \,,\ea$$
	i.e., $C$ is the solution of the discrete-time Lyapunov equation 
	\be C=ACA'+\Sigma\,. \label{eqn:covariance lyap eq} \ee 
	Consider now the following depomposition $A = DSD^{-1}$, where $D$ is a diagonal matrix s.t. $D:=[W\vect 1]_{\mc R \mc R}^{-1}$ and $S:=W_{\mc R \mc R}$ is a symmetric matrix.  Substituting  in \eqref{eqn:covariance lyap eq} we  obtain
	$$	C=DSD^{-1} CD^{-1}SD+\Sigma \,. $$
	Properly collecting terms, we get
	\[
	\begin{array}{rc} 
			&D(D^{-1}CD^{-1}-SD^{-1}CD^{-1}S)D = \Sigma \\
			&(I-S^2)D^{-1}CD^{-1}=D^{-1}\Sigma D^{-1} 
		\end{array} \]
	and consequently, by some extra computations,  $$C = \Sigma D^{-1}(I-S^2)^{-1}D\,.$$ 
	Substituting $S \!\! =\!\! D^{-1}AD$ and  $I \!\!=\!\! D^{-1}ID$  in the above, then
	\[
	\ba{rcl}C
	&=& \Sigma D^{-1}(D^{-1}ID-D^{-1}A^2D)^{-1}D \\
	&=&\Sigma D^{-1}(D^{-1}(I-A^2)D)^{-1}D\,. \ea \]
	The thesis follows.
\end{proof}
\begin{proof}[Proof of Corollary \ref{lemma: C properties}]
	
	(i) Given \eqref{eqn:covariance lyap eq}, we retrieve that, for any $x \in \mathbb{R}^{|\mc R|}\backslash \{ \vect 0\}$ it holds
	$$ x'Cx = (A'x)'C(A'x)+x'\Sigma x\,. $$
	Here, $(A'x)'C(A'x) \geq 0$, since $A \geq 0$ by definition and $C \geq 0$ as it is the limit  of a series of non negative terms. In addition, $x' \Sigma x > 0$ since $\Sigma$ it the covariance matrix of $V(t)$ whose variables $V_i(t)$ are independent by definition (see Chap. 7 in \cite{Matrix} for reference). Consequently, $x'Cx > 0$ for any $x \in \mathbb{R}^{|\mc R|}\backslash \{ \vect 0\}$ and thus, $C\succ 0$ by definition.
	\\(ii) It immediately follows from Proposition \ref{prop:general covariance}. As the precision matrix is defined as the inverse of the covariance matrix, it holds
	\be\label{eq:precision-matrix}H=C^{-1} = (I-A^2)\Sigma^{-1}\,,\ee
	where all  off-diagonal entries of $H$ are non-positive since $A\geq 0$.
\end{proof}
\subsection{Proof of Proposition \ref{prop:F formulation}} \label{app:Proposition 2}
Firstly, let us introduce a preliminary lemma that will be useful in the following.
\begin{lemma}
\label{lemma_proj}
Given $\mc K \subseteq \mc R$, it holds
\be \label{eqn:conditional expected}
\mathbb{E}[Y|X_{\mc K}] = \hat{\alpha}'X
\ee
where
$$	\hat{\alpha}_\mc K=(C_{\mc K \mc K})^{-1} (C \vect 1/n)_{\mc K} \quad ; \quad \hat{\alpha}_{-\mc K}=0\,. $$
\end{lemma}
\begin{proof}
From Result 7.7 in \cite{Mult_statistical_anal} it is well known that $	\mathbb{E}[Y|X_{\mc K}] = \hat{\alpha}'X$ with 
\be \label{eqn:alpha_min}	\hat{\alpha}=\underset{\beta: \mathrm{supp}(\beta)\subseteq\mc K}{\arg \min} \mathbb{E}[(Y-\beta'X)^2]\,.\ee
First, note that given $\mathrm{supp}(\hat \alpha)\subseteq \mc K$ then $ \hat{\alpha}_{-\mc K}=0$. Now, let us focus on $\hat{\alpha}_{\mc K}$. Using \eqref{eqn:alpha_min}, the following series of equations holds true:
\[
\begin{array}{ll}
\hat{\alpha}&=\underset{\beta: \mathrm{supp}(\beta)\subseteq\mc K}{\arg \min} \mathbb{E}[(Y-\beta'X)^2]\\
&=\underset{\beta: \mathrm{supp}(\beta)\subseteq\mc K}{\arg \min} \mathbb{E}[(((\vect 1/n)'-\beta')X)^2]\\
&=\underset{\beta: \mathrm{supp}(\beta)\subseteq\mc K}{\arg \min} \mathbb{E}[(\vect 1/n-\beta)'XX'(\vect 1/n-\beta)]\\
&=\underset{\beta: \mathrm{supp}(\beta)\subseteq\mc K}{\arg \min} (\vect 1/n-\beta)'C(\vect 1/n-\beta)\\
&=\underset{\beta: \mathrm{supp}(\beta)\subseteq\mc K}{\arg \min} -2\beta'C\vect 1/n+\beta'C\beta
\end{array}
\]
The minimum can so be found imposing that for any $k \in \mc K$
$$\frac{\partial (-2\hat{\alpha}'C\vect 1/n + \hat{\alpha}' C \hat{\alpha})}{\partial \hat{\alpha}_k} = 0$$
which is equivalent to 
$$(C\vect 1)_{\mc K} - (nC\hat{\alpha})_{\mc K} = 0 $$
and thus
$$C_{\mc K \mc K}(\vect 1_{\mc K}-n\hat{\alpha}_{\mc K})+C_{\mc K -\mc K}\vect 1_{-\mc K}=0\,.$$
Given that, from Corollary \ref{lemma: C properties}-i, $C \succ 0$, then  from Sylvester's criterion (see Theorem 7.2.5 in \cite{Matrix} for reference), $C_{\mc K \mc K}$ is invertible. The thesis follows.

\end{proof}
\begin{proof}[Proof of Proposition \ref{prop:F formulation}]
(i)	From law of total variances it holds that $ \Var(\mathbb{E}[Y|X_{\mc K}]) = \Var(Y) -  \mathbb{E}[\Var(Y|X_{\mc K})]$.
Then the result follows from Lemma \ref{lemma_proj}:
\[
\begin{split}
F(\mc K)&=\Var(\mathbb{E}[Y|X_{\mc K}])=\hat{\alpha}C\hat{\alpha} = \hat{\alpha}'_{\mc K}(C \hat{\alpha})_{\mc K}\\
&=\hat{\alpha}'_{\mc K}(C \vect 1/n)_{\mc K}= (C\vect 1/n)_{\mc K}'(C_{\mc K \mc K})^{-1}(C\vect 1/n)_{\mc K} \,.
\end{split}
\]
The thesis follows avoising constant terms.
\\(ii) First, proven \eqref{function_general}, we can write an extended expression for $F(\mc K)$:
\[
\begin{array}{rcl}	F(\mc K) =\!\!\!\!\!\!\!\!\!\!&&\vect 1_{\mc K}'C_{\mc K \mc K}\vect 1_{\mc K}+ \vect 1'_{\mc K}C_{\mc K -\mc K}\vect 1_{-\mc K}+\\[3pt]
&&\vect 1_{-\mc K}'C_{-\mc K \mc K} \vect 1_{\mc K} + \vect 1_{-\mc K}'C_{-\mc K \mc K} (C_{\mc K \mc K})^{-1}C_{\mc K - \mc K}\vect 1_{-\mc K} \,.
\end{array}
\]  
Thus, using \eqref{eqn:MSE-VarRed}, we retrieve that, avoiding constant terms,
\[
\begin{array}{rcl}
G(\mc K) \!\!\!\!\!\!\!&=&\!\!\!\! Var(Y)-F(\mc K)  \\
&=&\!\!\!\!\vect 1_{-\mc K}'(C_{-\mc K -\mc K} \!-\! C_{-\mc K \mc K}(C_{\mc K \mc K})^{-1} C_{\mc K - \mc K}) \vect 1_{-\mc K} \,.
\end{array}
\]
Based on theory about inverse of a partitioned nonsingular matrix (see [0.7.3.1] in \cite{Matrix}), it is well known that
\[
(H_{-\mc K-\mc K})^{-1}=C_{-\mc K -\mc K} - C_{-\mc K \mc K}(C_{\mc K \mc K})^{-1} C_{\mc K - \mc K}\,.
\]
The thesis follows.
\end{proof}
\subsection{Proof of Proposition \ref{prop:iterative step}-\ref{cor:iterative_inverse}}\label{app:Greedy}
\begin{proof}[Proof of Proposition \ref{prop:iterative step}]
	Recall that $$F({\mc K\cup \{i\}}) = (C \alpha)'_{\mc K\cup \{i\}}(C_{{\mc K\cup \{i\}} {\mc K\cup \{i\}}})^{-1}( C {\alpha})_{\mc K\cup \{i\}}$$	
	\textbf{Step 1.} \\ Note that, using the Schur Complement, it holds:
	$$(C_{{\mc K\cup \{i\}} {\mc K\cup \{i\}}})^{-1} = \begin{bmatrix}
			C_{\mc K \mc K} & C_{\mc K i} \\ C_{i \mc K} & C_{ii}
		\end{bmatrix}^{-1} = \begin{bmatrix}
			K_{11} & K_{12}\\ K_{21} & K_{22}
		\end{bmatrix}$$
	where \[
	\begin{array}{ll}
			K_{11}&=\left( C_{\mc K \mc K} - \frac{C_{\mc K i}C_{i \mc K}}{C_{ii}}\right)^{-1} \\
			K_{12}&= - K_{11} \frac{C_{\mc K i}}{C_{ii}}\\
			K_{21} &= - \frac{C_{\mc K i}}{C_{ii}} K_{11}\\
			K_{22} &= \frac{1}{C_{ii}} + \frac{C_{i \mc K}}{C_{ii}^2}K_{11}C_{\mc K i}
		\end{array}
	\]
	\textbf{Step 2. } \\ Note that, using Sherman-Morrison formula, it holds:
	\[
	\begin{array}{ll}
			K_{11}&=\left( C_{\mc K \mc K} - \frac{C_{\mc K i}C_{i \mc K}}{C_{ii}}\right)^{-1} =C_{\mc K \mc K}^{-1} + \frac{C_{\mc K \mc K}^{-1}C_{\mc K i}C_{i \mc K} C_{\mc K \mc K}^{-1}}{C_{ii}-C_{i \mc K}C_{\mc K \mc K}^{-1}C_{\mc K i}}
		\end{array}
	\]
	\textbf{Step 3.} \\
	We thus obtain:
	\[
	\begin{array}{ll}
			F({\mc K\cup \{i\}} )=&(C\alpha)_{\mc K}' (C_{\mc K \mc K})^{-1} (C\alpha)_{\mc K} \\
			&- (C\alpha)_{i}\frac{C_{\mc K i}}{C_{ii}}(C_{\mc K \mc K})^{-1}(C \alpha)_{\mc K} \\
			&- (C \alpha)_{\mc K}'(C_{\mc K \mc K})^{-1}\frac{C_{\mc K i}}{C_{ii}} (C\alpha)_{i}\\
			&+ \frac{(C\alpha)_i^2}{C_{ii}} + \frac{(C\alpha)_i^2}{C_{ii}}  C_{i \mc K}(C_{\mc K \mc K})^{-1}C_{\mc K i}
		\end{array}	
	\]
	The thesis follows.
\end{proof}
\begin{proof}[Proof of Proposition \ref{cor:iterative_inverse}]
	We apply here the same strategy proposed in Step 1-2 of proof of Proposition \ref{prop:iterative step}. No matrix inversion is required since $(C_{\mc K \mc K})^{-1}$ can be written as function of the inverse matrix at previous step. In addition, the first inverse matrix of the itarative procedure is computed for a single node observation, hence $(C_{ii})^{-1}=\frac{1}{C_{ii}}$.
\end{proof}
\end{document}